\documentclass[10pt]{amsart}
\usepackage{amsmath}
\usepackage{amsthm}
\usepackage[all]{xy}
\usepackage{verbatim}

\theoremstyle{plain}
\newtheorem{thm}{Theorem}

\newtheorem{lem}[thm]{Lemma}
\newtheorem{cor}[thm]{Corollary}

\newtheorem*{ithm}{Theorem}
\theoremstyle{remark}

\DeclareMathOperator{\aut}{Aut}

\DeclareMathOperator{\der}{Der}
\DeclareMathOperator{\enm}{End}

\DeclareMathOperator{\ic}{ic}
\DeclareMathOperator{\rc}{RC}
\DeclareMathOperator{\rh}{RH}
\DeclareMathOperator{\hmm}{Hom}

\newcommand{\n}{\mathbf{N_0}}

\newcommand{\z}{\mathbf{Z}}

\begin{document}

\title{Equivariant Lie-Rinehart cohomology}
\author{Eivind Eriksen}
\address{Oslo University College}
\email{eeriksen@hio.no}
\author{Trond St{\o}len Gustavsen}
\address{BI Norwegian School of Management}
\email{trond.s.gustavsen@bi.no}
\date{\today}

\begin{abstract}
In this paper, we study Lie-Rinehart cohomology for quotients of
singularities by finite groups, and interpret these cohomology groups
in terms of integrable connection on modules.
\end{abstract}

\maketitle

\section*{Introduction}

In Eriksen, Gustavsen \cite{eri-gus08}, we studied Lie-Rinehart
cohomology of singularities, and we gave an interpretation of these
cohomology groups in terms of integrable connections on modules of rank
one defined on the given singularities. The purpose of this paper is to
study equivariant Lie-Rinehart cohomology of the quotient of a
singularity by a finite group, and relate these cohomology groups to
integrable connections on modules of rank one defined on the quotient.

Let $k$ be an algebraically closed field of characteristic $0$, let $A$
be a reduced Noetherian $k$-algebra, and let $(M,\nabla)$ be a finitely
generated torsion free $A$-module of rank one with a (not necessarily
integrable) connection. There is an obstruction $\ic(M) \in
\rh^2(\der_k(A),\overline A)$ for the existence of an integrable
connection on $M$, where $\overline A = \enm_A(M)$, and if $\ic(M) =
0$, then $\rh^1(\der_k(A), \overline A)$ is a moduli space for the set
of integrable connections on $M$, up to analytical equivalence; see
Theorem \ref{t:eri-gus} in Eriksen, Gustavsen \cite{eri-gus08}.

We recall that an extension $R \subseteq S$ of a normal domain $R$ of
dimension two is called Galois if $S$ is the integral closure of $R$ in
$L$, where $K \subseteq L$ is a finite Galois extension of the quotient
field $K$ of $R$, and $R \subseteq S$ is unramified at all prime ideals
of height one.

\begin{ithm}
Let $G$ be a finite group, let $A$ be a normal domain of dimension two
of essential finite type over $k$ with a group action of $G$, and let
$M$ be an $A$-$G$ module that is finitely generated, maximal
Cohen-Macaulay of rank one as an $A$-module. We assume that the $M$
admits a connection. If $A^G \subseteq A$ is a Galois extension, then
the following hold:
\begin{enumerate}
\item We have that $\rh^n(\der_k(A^G),A^G) \cong
    \rh^n(\der_k(A),A)^G$ for all $n \ge 0$.
\item There is a canonical class $\ic(M^G) \in \rh^2(\der_k(A),
    A)^G$, called the integrability class, such that $\ic(M^G) = 0$
    if and only if there exists an integrable connection on $M^G$.
\item If $\ic(M^G)=0$, then $\rh^1(\der_k(A),A)^G$ is a moduli
    space for the set of integrable connections on $M^G$, up to
    analytical equivalence.
\end{enumerate}
\end{ithm}

We use this result to compute $\rh^{\ast}(\der_k(A^G),A^G)$ when $A =
k[x_1,x_2,x_3]/(f)$ is a quasi-homogeneous surface singularity and the
action of $G = \z_m$ on $A$ is of type $(m; m_1, m_2, m_3)$ such that
$A^G \subseteq A$ is a Galois extension. Explicitly, the action of $G$
on $A$ is given by
    \[ g{\ast}x_1 = \xi^{m_1} x_1, \quad g{\ast}x_2 = \xi^{m_2} x_2,
    \quad g{\ast}x_3 = \xi^{m_3} x_3 \]
for a cyclic generator $g$ of $G$ and a primitive $m$'th root of unity
$\xi \in k$. In the case when $f = x_1^3+x_2^3+x_3^3$ and $G = \z_3$
with action of type $(3; 1,1,2)$, we compute that
$\rh^n(\der_k(A^G),A^G) = 0$ for $n \ge 1$. In particular, if $M$ is an
maximal Cohen-Macaulay $A$-module of rank one that admits a connection,
then the $A^G$-module $M^G$ admits an integrable connection, unique up
to analytic equivalence.

\section{Basic definitions}

Let $k$ be an algebraically closed field of characteristic $0$, let $A$
be a commutative $k$-algebra, and let $M$ be an $A$-module. We define a
\emph{connection} on $M$ to be an $A$-linear map $\nabla: \der_k(A) \to
\enm_k(M)$ such that
\begin{equation*} \label{e:dp}
\nabla_D(am) = a \nabla_D(m) + D(a) \; m
\end{equation*}
for all $D \in \der_k(A), \; a \in A, \; m \in M$.

Let $\nabla$ be a connection on $M$. We define the \emph{curvature} of
$\nabla$ to be the $A$-linear map $R_{\nabla}: \der_k(A) \wedge
\der_k(A) \to \enm_A(M)$ given by
\begin{equation*}
R_{\nabla}(D \wedge D') = [ \nabla_D, \nabla_{D'} ] - \nabla_{[D,D']}
\end{equation*}
for all $D,D' \in \der_k(A)$. Notice that $R_{\nabla} \in
\hmm_A(\wedge^2_A \der_k(A), \enm_A(M))$. We say that $\nabla$ is an
\emph{integrable connection} if $R_{\nabla} = 0$, i.e. if $\nabla$ is a
homomorphism of Lie algebras.

For any integrable connection $\nabla$ on $M$, we consider the
Lie-Rinehart cohomology $\rh^{\ast}(\der_k(A),M,\nabla)$, see Eriksen,
Gustavsen \cite{eri-gus08}. We recall that the Lie-Rinehart cohomology
is the cohomology of the Lie-Rinehart complex, given by
    \[ \rc^n(\der_k(A),M,\nabla) = \hmm_A(\wedge_A^n \der_k(A),
    M) \]
for $n \ge 0$, with differentials $d^n: \rc^n(\der_k(A),M,\nabla) \to
\rc^{n+1}(\der_k(A),M,\nabla)$ given by
\begin{multline*}
    d^{n}(\xi)(D_{0} {\wedge} \cdots {\wedge} D_{n})
    = \sum_{0 \le i \le n} (-1)^{i} \;
    \nabla_{D_{i}}(\xi(D_{0} {\wedge} \cdots {\wedge}
    \widehat{D}_{i} {\wedge} \cdots {\wedge} D_{n})) \\
    + {\sum_{0\leq j<k\leq n} (-1)^{j+k} \xi([D_{j},D_{k}]
    {\wedge}D_{0} {\wedge} \cdots {\wedge} \widehat{D}_{j}
    {\wedge} \cdots {\wedge} \widehat{D}_{k}{\wedge} \cdots
    {\wedge}D_{n})}
\end{multline*}
for all $n \ge 0$ and all $\xi \in \rc^n(\der_k(A),M,\nabla), \; D_0,
D_1, \dots, D_n \in \der_k(A)$.

\section{Group actions}

Let $\sigma: G \to \aut_k(A)$ be a group action of a group $G$ on the
$k$-algebra $A$. For simplicity, we shall write $g{\ast}a =
\sigma(g)(a)$ for all $g \in G, \; a \in A$. We remark that $\sigma$
induces a group action of $G$ on $\der_k(A)$, given by
    \[ g{\ast}D = g D g^{-1} = \{ a \mapsto g{\ast}D(g^{-1}{\ast}a) \}
    \]
for $g \in G, D \in \der_k(A)$. Notice that we have $g{\ast}(aD) =
(g{\ast}a) \cdot (g{\ast}D)$ for all $a \in A$, $D \in \der_k(A)$.

We recall that an $A$-$G$ module structure on the $A$-module $M$ is a
group action $G \to \aut_k(M)$ that is compatible with the $A$-module
structure, i.e. a group action such that $g{\ast}(a m) = (g{\ast}a)
\cdot (g{\ast}m)$ for all $g \in G, \; a \in A, \; m \in M$. Hence
$\der_k(A)$ has a natural $A$-$G$ module structure induced by $\sigma$.

Let $M,N$ be $A$-$G$ modules, and consider the natural group actions of
$G$ on the $A$-modules $\hmm_A(M,N)$ and $M \wedge_A N$ given by
    \[ g{\ast}\phi = g \phi g^{-1} = \{ m \mapsto g{\ast}\phi(g^{-1}
   {\ast}m) \} \quad \text{and} \quad g{\ast}( m \wedge n ) = (g{\ast}m)
   \wedge (g{\ast}n) \]
for all $g \in G, \; \phi \in \hmm_A(M,N), \; m \in M$ and $n \in N$.
We remark that this gives $\hmm_A(M,N)$ and $M \wedge_A N$ natural
$A$-$G$ module structures.

Let $M$ be an $A$-$G$ module. Then there is an induced action of $G$ on
the set of connections on $M$, and $g{\ast}\nabla$ is given by
    \[ (g{\ast}\nabla)_D(m) = g{\ast}\nabla_{g^{-1}{\ast}D}(g^{-1}{\ast}m)
    \]
for any connection $\nabla: \der_k(A) \to \enm_k(M)$ and for all $g \in
G, \; D \in \der_k(A)$, $m \in M$. A straight-forward calculation shows
that
    \[ g{\ast}R_{\nabla} = R_{g{\ast}\nabla} \]
for all $g \in G$. In particular, if $\nabla$ is $G$-invariant, then
the same holds for $R_\nabla$.

If $G$ is a finite group, then there is a Reynolds' operator $M \to
M^G$ for any $A$-$G$ module $M$, given by
    \[ m \mapsto m' = \frac{1}{|G|} \, \sum_{g \in G} g{\ast}m \]
for all $m \in M$. Similarly, if $\nabla$ is a connection on an $A$-$G$
module $M$ and $G$ is a finite group, then
    \[ \nabla' = \frac{1}{|G|} \, \sum_{g \in G} g{\ast}\nabla \]
is a $G$-invariant connection on $M$. However, notice that $R_{\nabla_1
+ \nabla_2} \neq R_{\nabla_1} + R_{\nabla_2}$ in general. Hence, the
Reynolds' type operator $\nabla \mapsto \nabla'$ will not necessarily
preserve integrability.

Let $M$ be an $A$-$G$ module, and let $\nabla$ be an integrable
connection on $M$. Then $\rc^n(\der_k(A), M, \nabla) = \hmm_A(\wedge^n
\der_k(A), M)$ is an $A$-$G$ module for all $n \ge 0$. In fact, we have
the following result:

\begin{lem}
For any $g \in G$ such that $g{\ast}\nabla = \nabla$, the diagram
\[ \xymatrix{
M \ar[d]_{g{\ast}} \ar[r]^-{d^0} & \hmm_A(\der_k(A),
M) \ar[d]_{g{\ast}} \ar[r]^-{d^1} & \hmm_A(\wedge^2_A
\der_k(A), M) \ar[d]_{g{\ast}} \ar[r]^-{d^2} & \dots \\
M \ar[r]_-{d^0} & \hmm_A(\der_k(A), M) \ar[r]_-{d^1} &
\hmm_A(\wedge^2_A \der_k(A), M) \ar[r]_-{d^2} & \dots } \] commutes. In
particular, if $\nabla$ is $G$-invariant, then $G$ acts on the
Lie-Rinehart cohomology $\rh^{\ast}(\der_k(A),M,\nabla)$.
\end{lem}
\begin{proof}
Let $n \ge 0$ and let $\xi \in \rc^n(\der_k(A),M,\nabla)$. We must show
that if $g{\ast}\nabla = \nabla$, then $g{\ast}d^n(\xi) =
d^n(g{\ast}\xi)$, i.e. that $(g{\ast}d^n(\xi))(D_0 \wedge \dots \wedge
D_n) = (d^n(g{\ast}\xi))(D_0 \wedge \dots \wedge D_n)$ for all $D_0,
\dots, D_n \in \der_k(A)$. To simplify notation, we write
\begin{align*}
m_i &= g{\ast}\xi(g^{-1}{\ast}(D_0 \wedge \dots \widehat{D_i} \wedge
\dots \wedge D_n)) \\
m_{jk} &= g{\ast}\xi(g^{-1}{\ast}([ D_j, D_k] \wedge D_0 \wedge
\dots \wedge \widehat{D_j} \wedge \dots \wedge \widehat{D_k} \wedge
\dots \wedge D_n))
\end{align*}
for $0 \le i \le n$ and $0 \le j < k \le n$. Using this notation, we
compute that
\begin{align*}
(g{\ast}d^n(\xi))(D_0 \wedge \dots \wedge D_n)
&= g{\ast}d^n(\xi)(g^{-1}{\ast}(D_0 \wedge \dots \wedge D_n)) \\
&= \sum_{i} (-1)^i \; g{\ast}\nabla_{g^{-1}{\ast}D_i}(g^{-1}{\ast}m_i)
+ \sum_{j<k} (-1)^{j+k} \; m_{jk} \\
&= \sum_{i} (-1)^i \; (g{\ast}\nabla)_{D_i}(m_i) + \sum_{j<k}
(-1)^{j+k} \; m_{jk}
\end{align*}
and that
\begin{align*}
& (d^n(g{\ast}\xi))(D_0 \wedge \dots \wedge D_n) = \sum_{i} (-1)^i \;
\nabla_{D_i}((g{\ast}\xi)(D_0 \wedge \dots \widehat{D_i} \wedge \dots
D_n)) \\
&+ \sum_{j<k} (-1)^{j+k} \; (g{\ast}\xi)([D_j, D_k] \wedge D_0 \wedge
\dots \wedge \widehat{D_j} \wedge \dots \wedge \widehat{D_k} \wedge
\dots \wedge D_n) & \\
&= \sum_i (-1)^i \; \nabla_{D_i}(m_i) + \sum_{j<k} (-1)^{j+k} \;
m_{jk} &
\end{align*}
since $g^{-1}[D_j,D_k] = [g^{-1}D_j, g^{-1}D_k]$. But $g{\ast}\nabla =
\nabla$ by assumption, and it follows that $g{\ast}d^n(\xi) =
d^n(g{\ast}\xi)$.
\end{proof}

\section{Integrable connections and Lie-Rinehart cohomology}

In the rest of this paper, we assume that $A$ is a reduced Noetherian
$k$-algebra and that $M$ is a finitely generated torsion free
$A$-module of rank one. Hence there is an isomorphism $f: Q(A) \to Q(A)
\otimes_A M$ of $Q(A)$-modules, where $Q(A)$ is the total ring of
fractions of $A$. We identify $M$ with its image in $Q(A) \otimes_A M$,
and let $N = f^{-1}(M) \subseteq Q(A)$, so that $f: N \to M$ is an
isomorphism of $A$-modules. Then there is an identification
    \[ \enm_A(M) \cong \overline A = \{ q \in Q(A): q \cdot N
    \subseteq N \} \]
We consider $\overline A$ as a commutative overring with $A \subseteq
\overline A \subseteq Q(A)$. If $A$ is normal, then $\overline A = A$,
see Eriksen, Gustavsen \cite{eri-gus08}, Proposition 3.1.

We remark that if $\nabla$ is a connection on $M$, then there is an
induced connection $\overline \nabla$ on the $A$-module $\overline A
\cong \enm_A(M)$. The induced connection $\overline \nabla$ is trivial
in the sense that
    \[ \overline \nabla_D(q) = \overline D(q) \]
for any $D \in \der_k(A)$ and any $q \in \overline A$, where $\overline
D$ is the natural lifting of $D$ to $\der_k(Q(A))$. In particular,
$\overline \nabla$ is an integrable connection on $\overline A$.

Since $\overline \nabla$ is a canonical integrable connection on
$\overline A$, it is natural to consider the Lie-Rinehart cohomology
$\rh^{\ast}(\der_k(A),\overline A)$ with values in $(\overline A,
\overline \nabla)$. We recall the following result:

\begin{thm}[Eriksen-Gustavsen] \label{t:eri-gus}
Let $A$ be a reduced Noetherian $k$-algebra and let $M$ be a finitely
generated, torsion free $A$-module of rank one. We assume that $M$
admits a connection.
\begin{enumerate}
\item There is a canonical class $\ic(M) \in \rh^2(\der_k(A),
    \overline A)$, called the integrability class, such that
    $\ic(M) = 0$ if and only if there exists an integrable
    connection on $M$.
\item If $\ic(M)=0$, then $\rh^1(\der_k(A),\overline A)$ is a
    moduli space for the set of integrable connections on $M$, up
    to analytical equivalence.
\end{enumerate}
\end{thm}
\begin{proof}
See Eriksen, Gustavsen \cite{eri-gus08}, Proposition 3.2 and Theorem
3.4.
\end{proof}

\section{The equivariant case}

Let $A$ be a reduced Noetherian $k$-algebra with a group action
$\sigma: G \to \aut_k(A)$, and let $M$ be an $A$-$G$ module that is
finitely generated, torsion free of rank one as an $A$-module. Then
there is a natural group action of $G$ on $Q(A)$ induced by $\sigma$,
given by
    \[ g{\ast}q = g{\ast}\left( \frac{a}{s} \right) =
    \frac{g{\ast}a}{g{\ast}s} \]
for any $g \in G$ and any $q = a/s \in Q(A)$. Let us consider the (not
necessarily equivariant) isomorphism $f: Q(A) \to Q(A) \otimes_A M$ of
$Q(A)$-modules, and the submodule $N = f^{-1}(M) \subseteq Q(A)$. Then
we have an identification
    \[ \enm_A(M) \cong \overline A = \{ q \in Q(A): q \cdot N \subseteq
    N \} \]
as above, and the group action of $G$ on $\enm_A(M)$ induced by the
$A$-$G$ module structure on $M$ coincides with the natural group action
on $\overline A$ induced by the group action of $G$ on $Q(A)$. In fact,
if $\phi \in \enm_A(M)$ corresponds to $q \in \overline A$, then $\phi$
is given by $\phi(m) = q \cdot m$ for all $m \in M$, and
    \[ (g{\ast}\phi)(m) = g{\ast}\phi(g^{-1}{\ast}m) = g{\ast}(q \cdot
    (g^{-1}{\ast}m)) = (g{\ast}q) \cdot m \]
for all $g \in G$ and $m \in M$. We also remark that since the induced
connection $\overline \nabla$ on $\overline A$ is trivial, it is clear
that $\overline \nabla$ is $G$-invariant.

Assume that $R$ is a normal domain of dimension two over $k$, let $K
\subseteq L$ be a finite Galois extension of the quotient field $K$ of
$R$, and let $S$ be the integral closure of $R$ in $L$. If the
extension $R \subseteq S$ is unramified at all prime ideals of height
one, we say that it is a \emph{Galois extension}. It is known that if
$S = k[x,y]$ is a polynomial ring and $G \subseteq \aut_k(S)$ is a
finite subgroup without pseudo-reflections, then $S^G \subseteq S$ is a
Galois extension.

\begin{thm}
Let $G$ be a finite group, let $A$ be a normal domain of dimension two
of essential finite type over $k$ with a group action of $G$, and let
$M$ be an $A$-$G$ module that is finitely generated, maximal
Cohen-Macaulay of rank one as an $A$-module. We assume that the $M$
admits a connection. If $A^G \subseteq A$ is a Galois extension, then
the following hold:
\begin{enumerate}
\item We have that $\rh^n(\der_k(A^G),A^G) \cong
    \rh^n(\der_k(A),A)^G$ for all $n \ge 0$.
\item There is a canonical class $\ic(M^G) \in \rh^2(\der_k(A),
    A)^G$, called the integrability class, such that $\ic(M^G) = 0$
    if and only if there exists an integrable connection on $M^G$.
\item If $\ic(M^G)=0$, then $\rh^1(\der_k(A),A)^G$ is a moduli
    space for the set of integrable connections on $M^G$, up to
    analytical equivalence.
\end{enumerate}
\end{thm}
\begin{proof}
Since $A$ is normal, we have that $\overline A = A$, and therefore
$\rh^{\ast}(\der_k(A),A)$ is the Lie-Rinehart cohomology associated
with $(\overline A, \overline \nabla)$. The functor $M \mapsto M^G$ is
exact since $G$ is a finite group, hence $\rh^{\ast}(\der_k(A),A)^G
\cong \text{H}^{\ast}(\rc^{\ast}( \der_k(A),A)^G)$. Moreover, it
follows as in Proposition 4.4 of Gustavsen, Ile \cite{gus-ile08} that
    $$\hmm_A(\wedge^n \der_k(A),A)^G \cong \hmm_{A^G}(\wedge^n
    \der_k(A^G),A^G)$$
and this proves the first part of the theorem. By assumption, $A^G$ is
a normal domain, and we notice that $M^G$ is an maximal Cohen-Macaulay
$A^G$-module, see Gustavsen, Ile \cite{gus-ile08}, Proposition 4.3.
Therefore
    $$\nabla' = \frac{1}{|G|} \; \sum_{g \in G} \, g{\ast}\nabla$$
is a connection on $M^G$, and $\rh^n(\der_k(A^G),A^G)$ is the
Lie-Rinehart cohomology associated with $(\overline{A^G},
\overline{\nabla'})$. The rest of the theorem follows from Theorem
\ref{t:eri-gus}.
\end{proof}

We remark that this result can be generalized to higher dimensions with
suitable conditions on the extension $A^G \subseteq A$. In the case of
surface quotient singularities, we have the following corollary:

\begin{cor}
Let $A = k[x_1, x_2]$ be a polynomial algebra, and let $G \subseteq
\aut_k(A)$ be a finite subgroup without pseudo-reflections. Then we
have $\rh^n(\der_k(A^G),A^G) = 0$ for $n \ge 1$. In particular, any
maximal Cohen-Macaulay module over $A^G$ of rank one has an integrable
connection, unique up to analytic equivalence.
\end{cor}

\section{Quotients of quasi-homogeneous singularities}

Let $A = k[x_1,x_2,x_3]/(f)$ be an integral quasi-homogeneous surface
singularity. We write
    \[ f = \sum_{\alpha \in \n^3} \; c_{\alpha} \,
    {\underline x}^{\alpha} \]
and define $I(f) = \{ \alpha \in \n^3: c_{\alpha} \neq 0 \}$. Then
there are integral weights $d = \deg(f)$ and $d_i = \deg(x_i)$ for $i =
1,2,3$ such that $\alpha_1 d_1 + \alpha_2 d_2 + \alpha_3 d_3 = d$ for
all $\alpha \in I(f)$.

Let $(m_1,m_2,m_3) \in \n^3$ and let $G = \z_m = \langle g: g^m = 1
\rangle$ be the cyclic group of order $m$. If $\alpha_1 m_1 + \alpha_2
m_2 + \alpha_3 m_3 = m$ for all $\alpha \in I(f)$, then there is a
group action of $G$ on $A$ given by
\begin{align*}
g{\ast}x_1 &= \xi^{m_1} \cdot x_1 \\
g{\ast}x_2 &= \xi^{m_2} \cdot x_2 \\
g{\ast}x_3 &= \xi^{m_3} \cdot x_3
\end{align*}
where $\xi \in k$ is a primitive $m$'th root of unity. We call this a
group action of $G$ on $A$ of type $(m;m_1,m_2,m_3)$.

\begin{thm} \label{t:qh-rch}
Let $A = k[x_1,x_2,x_3]/(f)$ be an integral quasi-homogeneous surface
singularity with weights $(d; d_1,d_2,d_3)$, and let $G$ be a cyclic
group of order $m$ with a group action on $A$ of type
$(m;m_1,m_2,m_3)$. If $A^G \subseteq A$ is a Galois extension, then we
have
\begin{enumerate}
    \item $\rh^0(\der_k(A),A)^G = \left( A_0 \right)^G = k$
    \item $\rh^1(\der_k(A),A)^G = \left( \rh^1(\der_k(A),A)_0
        \right)^G = \left( A_{d-d_1-d_2-d_3} \cdot e_1 \right)^G$
    \item $\rh^2(\der_k(A),A)^G = \left( \rh^2(\der_k(A),A)_0
        \right)^G = \left( A_{d-d_1-d_2-d_3} \cdot e_2 \right)^G$
\end{enumerate}
for a generator $e_n \in \rh^n(\der_k(A),A)_0$ with $g{\ast}e_n =
\xi^{m_1+m_2+m_3-m} \cdot e_n$ for $n = 1,2$.
\end{thm}
\begin{proof}
In Theorem 6.2 of Eriksen, Gustavsen \cite{eri-gus08}, we considered
the case when $A = k[x_1,x_2,x_3]/(f)$ is quasi-homogeneous, and proved
that
\begin{align*}
\rh^0(\der_k(A),A) &= A_0 = k \\
\rh^1(\der_k(A),A) &= \rh^1(\der_k(A),A)_0 = A_{d-d_1-d_2-d_3} \cdot
\psi^{(4)} \\
\rh^2(\der_k(A),A) &= \rh^2(\der_k(A),A)_0 = A_{d-d_1-d_2-d_3} \cdot
\Delta^*
\end{align*}
An explicit description of $e_1 = \psi^{(4)}$ and $e_2 = \Delta^*$ is
given in Eriksen, Gustavsen \cite{eri-gus08}, Section 6, and we see
that
    \[ g{\ast}\psi^{(4)} = \xi^{m_1+m_2+m_3-m} \cdot \psi^{(4)}, \quad
    g{\ast}\Delta^* = \xi^{m_1+m_2+m_3-m} \cdot \Delta^* \]
\end{proof}

\subsection*{Example}

Let $A = k[x_1,x_2,x_3]/(f)$ with $f = x_1^3+x_2^3+x_3^3$, and consider
the action of $G = \z_3$ on $A$ of type $(3;1,1,2)$, given by
$g{\ast}x_i = \xi \cdot x_i$ for $i=1,2$ and $g{\ast}x_3 = \xi^2 \cdot
x_3$. In this case, it is known that $A^G \subseteq A$ is a Galois
extension, see Gustavsen, Ile \cite{gus-ile08}. We have $d-d_1-d_2-d_3
= 0$ and $m_1+m_2+m_3-m = 1$, hence it follows from Theorem
\ref{t:qh-rch} that $G$ acts non-trivially on $\rh^n(\der_k(A),A) \cong
A_0 = k$ for $n = 1,2$. This implies that $\rh^n(\der_k(A^G),A^G) \cong
\rh^n(\der_k(A),A)^G = 0$ for $n = 1,2$. In particular, if $M$ is a
maximal Cohen-Macaulay $A$-module of rank one that admits a connection,
then $M^G$ admits an integrable connection, unique up to analytic
equivalence.

In fact, it is known that $A^G$ is a rational singularity in this case,
and it follows from Remark 5.2 in Eriksen, Gustavsen \cite{eri-gus08}
that $\rh^1(\der_k(A^G),A^G) = 0$. On the other hand, the result that
$\rh^2(\der_k(A^G),A^G) = 0$ is stronger than the results obtained in
\cite{eri-gus08}.

\bibliographystyle{amsplain}
\bibliography{eeriksen}

\end{document}